\documentclass{birkjour}

\usepackage{amssymb,amsfonts,amsmath,latexsym}
\usepackage{xfrac}

 \newtheorem{thm}{Theorem}[section]

 \theoremstyle{definition}
 
 \theoremstyle{remark}

 \numberwithin{equation}{section}
 \newtheorem{ax}{A}
 \newtheorem{lemma}{Lemma}
 
 \newcommand{\epr}{\hfill $\Box$\mbox{}\\ }

\newcommand{\onei}{\ensuremath{\mathbf{1}^{i}}}
\newcommand{\onej}{\ensuremath{\mathbf{1}^{j}}}

\newcommand{\onel}{\ensuremath{\mathbf{1}^{l}}}

\newcommand{\onen}{\ensuremath{\mathbf{1}^{n}}}
\newcommand{\oner}{\ensuremath{\mathbf{1}^{r}}}
\newcommand{\ones}{\ensuremath{\mathbf{1}^{s}}}

\newcommand{\oneone}{\ensuremath{\mathbf{1}^{1}}}
\newcommand{\onetwo}{\ensuremath{\mathbf{1}^{2}}}

\newcommand{\RR}{\ensuremath{\mathbb{R}}}
\newcommand{\QQ}{\ensuremath{\mathbb{Q}}}
\newcommand{\NN}{\ensuremath{\mathbb{N}}}

\newcommand{\Hom}{Homogeneity}
\newcommand{\AQ}{Deviations Balancedness}
\newcommand{\IE}{Inverse Effects}
\newcommand{\RIE}{Restricted Inverse Effects}

\providecommand{\abs}[1]{\left\lvert #1\right\rvert}

\DeclareMathOperator*{\Int}{int}

\begin{document}

%

\title{Axiomatic characterization of the $\chi^{2}$ dissimilarity measure}

\author{Denis Bouyssou}

\address{%
CNRS and Universit\'{e} Paris Dauphine}
\email{bouyssou@lamsade.dauphine.fr}

\thanks{Authors are listed alphabetically and have contributed equally.}

\author{Thierry Marchant}
\address{Ghent University}
\email{thierry.marchant@ugent.be}

\author{Marc Pirlot}
\address{UMons}
\email{marc.pirlot@umons.ac.be}

\subjclass{Primary 39B22; Secondary 62G10}

\keywords{Chi squared, dissimilarity measure}

\date{\today}


\begin{abstract}
We axiomatically  characterize the $\chi^{2}$ dissimilarity measure. To this end, we solve a new generalization of a functional equation discussed in Aczel (Lectures on functional equations and their applications,  Academic Press, 1966).
\end{abstract}

\maketitle

\section{Introduction}

Let $N= \{1, 2, \ldots, n\}$ be a set of categories (with $n \geq 2$).  The vector $x=(x_{1}, \ldots, x_{n})$ represents the respective numbers of observations in each category and the total number of observations is denoted by $s(x)= \sum_{i \in N} x_{i}$.
We want to measure the dissimilarity between the observed distribution $x$ and a reference distribution $\pi = (\pi_{1}, \ldots, \pi_{n})$, with $\sum_{i \in N} \pi_{i} = 1$ and $\pi_{i} \in \mathbb{Q}_{++}$ for all $i \in N$, where $\mathbb{Q}_{++}$ is the set of positive rational numbers. We exclude reference distributions with null components because the $\chi^{2}$ dissimilarity measure is not defined when a component is zero.
The set of all observed distributions is $X= \mathbb{N}^{N}$, i.e. the set of all mappings from $N$ to $\mathbb{N}_{0}$, where $\NN_{0}$ is the set of non-negative integers. The set $\Pi$ of all reference distributions is  defined by $\Pi = \{\pi \in \mathbb{Q}_{++}^{N} : \sum_{i \in N} \pi_{i} = 1 \}$. 

A dissimilarity measure $f$ is a mapping from $X \times \Pi$ to $\mathbb{R}_{+}$ (the set of non-negative real numbers) satisfying 
 $f(x, \pi)=0$ iff $x/s(x) = \pi$. 
 It measures how far  the observed distribution is from the reference.   In this paper, we axiomatically characterize the $\chi^{2}_{1}$ dissimilarity measure defined by
$$\chi^{2}_{1}(x, \pi) = \sum_{i \in N} \frac{(s(x) \pi_{i} - x_{i} )^{2}}{s(x) \pi_{i}}$$
and frequently used in statistics as a measure of goodness of fit. 

The dissimilarity measure $\chi^{2}_{0}$ defined by $\chi^{2}_{0}(x,\pi)= \chi^{2}_{1}(x,\pi)/s(x)$ has been characterized in \cite{KaufmanMathaiRathie1972} and we will also provide a new characterization thereof. It is popular  in ecology \cite{Greenacre2017}, sociology \cite{ReardonFirebaugh2002}, economics \cite{Yalonetzky2012}, and so on.

While we consider in our paper that the number $n$ of categories is given and fixed, \cite{KaufmanMathaiRathie1972} considers that $n$ can vary. Depending on the context, one or the other assumption can be more relevant. For instance, when we use  Pearson's $\chi^{2}$ test, we have a sample distributed over $n$ categories and the $p$-value is computed conditional on a theoretical probability distribution with the same number $n$ of categories. If we repeat the experiment and draw other samples, we obtain other $p$-values always based on the same theoretical probability distribution with the same number $n$ of categories. It therefore makes sense to consider $n$ as given.

A common feature of \cite{KaufmanMathaiRathie1972} and our paper is that 
we use a framework in which $\pi$ can vary and such that comparisons of the dissimilarity measure across different reference distributions are relevant. Yet, unlike \cite{KaufmanMathaiRathie1972}, we also consider the case in which the reference distribution $\pi$ is fixed (as in our Pearson's $\chi^{2}$ example).

For characterizations of other dissimilarity measures, in the context of political sciences, see \cite{BouyssouMarchantPirlotAOR2020}. See \cite{BertoluzzaBaccoDoldi2004} for a characterization of a wide class of dissimilarity measures.
While we  consider dissimilarity measures in this paper, it is also interesting to consider dissimilarity rankings as in \cite{CowellDavidsonFlachaire2015}. 

Section 2 presents our main conditions and results. Section 3 shows the independence of the conditions used in our results. Section 4 concludes the discussion. All the proofs are gathered in Section 5.

\section{Axioms and results}

The dissimilarity measures $\chi_{0}^{2}$ and $\chi_{1}^{2}$  are  homogeneous of degree 0 and 1, respectively, where homogeneity is defined as follows.%
\begin{ax}
\emph{\Hom\ of degree $\omega$.} 
\label{A:Hom}
For all positive integers $\lambda$ and $x \in X$,  $f(\lambda x,\pi) = \lambda^{\omega} f(x,\pi)$.
\end{ax}
In statistics, it seems unanimously accepted that a dissimilarity measure (used as a goodness-of-fit statistic) should be homogeneous of degree 1, but in ecology, many researchers seem to favour homogeneity of degree 0. Indeed, when they measure the dissimilarity between the species distribution in an ecosystem and a reference distribution, they want the dissimilarity to be independent of the size of the ecosystem. It is easy to see that Homogeneity of degree 0 (resp.\ 1) is satisfied by $\chi^2_0$ (resp.\ $\chi^2_1$). Indeed, we have
$$\chi^{2}_{0}(\lambda x, \pi) = \sum_{i \in N} \frac{( \pi_{i} - \lambda x_{i}/  s( \lambda x) )^{2}}{ \pi_{i}} = \sum_{i \in N} \frac{( \pi_{i} -  x_{i}/  s(x) )^{2}}{ \pi_{i}} = \chi^{2}_{0}( x, \pi)$$
and
$$\chi^{2}_{1}(\lambda x, \pi) = \sum_{i \in N} \frac{(  s(\lambda x) \pi_{i} - \lambda x_{i})^{2}}{ s(\lambda x) \pi_{i}} = \lambda \sum_{i \in N}  \frac{( s(x) \pi_{i} -  x_{i} )^{2}}{ s(x) \pi_{i}} = \lambda\chi^{2}_{1}( x, \pi).$$

Suppose the dissimilarity between a distribution $x$ and $\pi$ is zero. This implies $x = k \pi$ for some positive integer $k$.
The next condition  states that, when we modify $k \pi$ by moving a single individual from  category $l$ to $j$, then the dissimilarity measure is inversely proportional to the harmonic mean  of $\pi_{j}$ and $\pi_{l}$. 
Let $\onei \in X$ be a vector such that $\onei_{i}=1$ and $\onei_{j}=0$ for all $j \neq i$.
\begin{ax}
\emph{\IE.} 
\label{A:IE}
If $k \pi , k \pi' \in X$, then, for all $j,l,r,s \in N$, with $j \neq l$ and $r \neq s$,
$$\frac{f(k \pi + \onej -\onel, \pi)}{f(k \pi' + \oner - \ones, \pi')} = \frac{\frac{1}{\pi_{j}} + \frac{1}{\pi_{l}}}{ \frac{1}{\pi'_{r}} + \frac{1}{\pi'_{s}}}.$$
\end{ax}

 In our first result, we will use a restricted variant of \IE\ in which $\pi=\pi'$. This weaker condition is named \RIE\ and is trivially satisfied when $n=2$.  We now prove that \IE\ is satisfied by $\chi^{2}_{0}$:
$$\chi^{2}_{0}(k \pi + \onej -\onel, \pi) =  \frac{( -1/  k )^{2}}{ \pi_{j}} + \frac{( 1/  k )^{2}}{ \pi_{l}} = \frac{1}{k^{2}} \left(\frac{1}{\pi_{j}} + \frac{1}{\pi_{l}} \right).$$
The proof for $\chi^{2}_{1}$ is similar.

Let $x$ and $y$ be two observed distributions of size $k$. 
The deviation between  $x$ and  $k\pi$ is $x - k\pi$. The corresponding deviation for $y$ is $y - k\pi$. If we add these two vectors of deviations, we obtain $x+y-2k\pi$ and the corresponding observed distribution is $x+y-2k\pi+k\pi = x+y-k\pi$ (provided all components are non-negative). Hence,  $f(x+y-k \pi, \pi)$ represents the dissimilarity corresponding to the additive combination of two deviations: between $x$ (resp.\ $y$) and $k\pi$.
Similarly, $f(x-y+k \pi, \pi)$  corresponds to the subtractive combination of the same two deviations. Finally, $f(x+y-k \pi, \pi) + f(x-y+k \pi, \pi)$ corresponds in some sense to four deviations (two $x$- and two $y$-deviations) combined once additively and once subtractively. Our next condition states that this must be equal to  $2f(x, \pi) +  2f(y, \pi)$, which is another way to combine the same four deviations.
\begin{ax}
\emph{\AQ.} 
\label{A:Q}
For all $x,y \in X$ with $s(x) = s(y) = k$,  if $x+y-k \pi \in X$  and $x-y+k \pi \in X$, then 
$$f(x+y-k \pi, \pi) + f(x-y+k \pi, \pi) =  2 \big( f(x, \pi) +  f(y, \pi) \big).$$
\end{ax}
This  condition is  inspired by \cite{DAgostinoDardanoni2009}, in which they characterize the Euclidean distance in $\RR^{n}$.
Let us prove that $\chi^{2}_{1}$ satisfies \AQ. We have 
\begin{align*}
\chi^{2}_{1}(x+y-k \pi, \pi) & = \sum_{i \in N} \frac{\big(s(x+y-k \pi) \pi_{i} - (x_{i}+y_{i}-k \pi_{i} ) \big)^{2}}{s(x+y-k \pi) \pi_{i}} \\
& = \sum_{i \in N} \frac{(2k \pi_{i} - x_{i}-y_{i} )^{2}}{k \pi_{i}} 
\end{align*}
and 
\begin{align*}
\chi^{2}_{1}(x-y+k \pi, \pi) & = \sum_{i \in N} \frac{\big(s(x-y+k \pi) \pi_{i} - (x_{i}-y_{i}+k \pi_{i} ) \big)^{2}}{s(x-y+k \pi) \pi_{i}} \\
& = \sum_{i \in N} \frac{( x_{i}-y_{i} )^{2}}{k \pi_{i}} .
\end{align*}
Hence, $\chi^{2}_{1}(x+y-k \pi, \pi)+\chi^{2}_{1}(x-y+k \pi, \pi)$ is equal to 
\begin{align*}
 & \sum_{i \in N} \frac{(2k \pi_{i} - x_{i}-y_{i} )^{2}}{k \pi_{i}} + \sum_{i \in N} \frac{( x_{i}-y_{i} )^{2}}{k \pi_{i}} \\
 = &   \sum_{i \in N} \frac{2(k^{2} \pi^{2}_{i} + x_{i}^{2} - 2k \pi_{i} x_{i}) + 2(k^{2} \pi^{2}_{i} + y_{i}^{2} - 2k \pi_{i} y_{i}) }{k \pi_{i}} \\
  = & 2 \chi^{2}_{1}(x, \pi) + 2 \chi^{2}_{1}(y, \pi).
 \end{align*}

We are now ready to state our first result in which we  consider that  $\pi$ is given and does not vary.

\begin{thm}
\label{Chi2}
Assume $\pi$ is given.  For $\omega \in \{0,1\}$,
a dissimilarity measure $f$ satisfies \Hom\ of degree $\omega$,  \AQ\ and \RIE\ iff 
$f = \gamma  \chi^{2}_{\omega}$,
for some positive $\gamma \in \mathbb{R}$.  \RIE\ is not required when $n=2$.
\end{thm}

Notice that Theorem~\ref{Chi2} does not hold when $\pi$ is not fixed. Indeed, for any $\phi: \Pi \to \RR_{+}$ with $\phi$ not constant, the dissimilarity measure 
$$f_{\phi}(x,\pi)= \phi(\pi) \sum_{i \in N} \frac{( \pi_{i} - x_{i}/s(x) )^{2}}{ \pi_{i}}
\label{f.phi}$$
satisfies \Hom\ of degree 1, \AQ\ and \RIE\  but is not of the form $\gamma \chi^{2}_{0}$ or $\gamma \chi^{2}_{1}$.
In order to  characterize  the dissimilarity measure $\chi^{2}$  when $\pi$  varies, we need the full power of \IE.

\begin{thm}
\label{Chi2.pi.variable}
For $\omega \in \{0,1\}$,
a dissimilarity measure $f$ satisfies \Hom\ of degree $\omega$,  \AQ\ and \IE\ iff 
$f = \gamma  \chi^{2}_{\omega}$,
for some positive $\gamma \in \mathbb{R}$.  
\end{thm}

\section{Independence of the axioms}

In order to prove the independence of the conditions characterizing $\chi^{2}_{0}$ with variable $\pi$ , we provide three examples of dissimilarity measures violating only one of the three conditions  in Theorem~\ref{Chi2.pi.variable}. 

The dissimilarity measure $\chi^{2}_{1}$  violates  \Hom\ of degree 0 but satisfies \AQ\ and \IE. The dissimilarity measure
$$f(x,\pi)  = \sum_{i \in N} \frac{\abs{ \pi_{i} - x_{i}/s(x) }}{ \pi_{i} }$$
 violates  \AQ\ but satisfies \Hom\ of degree 0 and \IE. 
 The dissimilarity measure
$$f(x,\pi)  =  \sum_{i \in N} ( \pi_{i} - x_{i} / s(x) )^{2}$$
violates \IE\ but satisfies \Hom\ of degree 0 and \AQ. 

Our examples are easily adapted to prove the independence of the conditions characterizing $\chi^{2}_{1}$ with variable $\pi$ . Finally, our examples can also be used for Theorem~\ref{Chi2} since it involves the same conditions as Theorem~\ref{Chi2.pi.variable} except for \RIE\ which is weaker than \IE.

\section{Discussion}
\label{sec.discussion}

Theorems~\ref{Chi2} and~\ref{Chi2.pi.variable} characterize the dissimilarity measures $\chi_{0}^{2}$ and $\chi_{1}^{2}$ up to a multiplication by a positive real number $\gamma$. We could easily add a condition characterizing exactly $\chi_{0}^{2}$ or $\chi_{1}^{2}$. For instance, the extra condition $f\big( \oneone,(1/n,  \ldots, 1/n) \big) = n-1$ is enough to force $\gamma=1$ in both  characterizations. Yet, unlike
\cite{KaufmanMathaiRathie1972}, we consider that such a normalization is not really interesting. Indeed $\chi_{1}^{2}$ and $\gamma \chi_{1}^{2}$ (with $\gamma \neq 1$) convey exactly the same information, just like a distance measurement in meters or yards. In particular, if we want to  perform a Pearson's $\chi^{2}$ test, we are free to use Pearson's statistic (i.e.\ $\chi_{1}^{2}$) and to compute the $p$-value using the $\chi^{2}$ density or to use $\gamma \chi_{1}^{2}$ (with an arbitrary $\gamma$) and to compute the $p$-value using the corresponding  density. The resulting $p$-value will of course be identical. 
The same holds for $\chi_{0}^{2}$ and $\gamma \chi_{0}^{2}$.

\section{Proofs}

\bigskip
We need a few lemmas before proving Theorem~\ref{Chi2}.

\begin{lemma}
\label{lemma:0-1-Homogeneity}
Let $f(x,\pi) = s(x) f'(x,\pi)$. Then 
 $f$ satisfies \Hom\ of degree 1  iff $f'$ satisfies \Hom\ of degree 0. And $f$ satisifies \AQ\ (resp.\ \IE) iff $f'$ satisfies \AQ\ (resp.\  \IE).
\end{lemma}

\begin{proof}
Since $f$ satisfies \Hom\ of degree 1, we have   $f(\lambda x,\pi)= \lambda f(x,\pi)$ for all positive integers $\lambda$. We thus have $\lambda s(x) f'(\lambda x, \pi) = \lambda s(x) f'(x,\pi)$. Hence $f'(\lambda x, \pi) = f'(x,\pi)$ and $f'$ is homogeneous of degree 0. The proof of the reverse implication is similar. 
The rest of the proof is left to the reader.
\end{proof}

\begin{lemma}
\label{lemma:0-Homogeneity}
Suppose $\pi$ is fixed.
If a dissimilarity measure $f$ satisfies \Hom\ of degree 0,  then $f(x, \pi)= F(x/s(x))$, for some mapping $F: \Pi \rightarrow \mathbb{R}_{+}$.
\end{lemma}

\begin{proof}
Since $\pi$ is fixed, we can define a mapping $g: X \rightarrow \mathbb{R}_{+}$ such that $f(x, \pi)= g(x)$.  Define now the mapping $F: \Pi \rightarrow \mathbb{R}_{+}$ as follows.  For any $p \in \Pi$, $F(p) = g(x)$ if there is $x \in X$ such that $p = x/s(x)$. The mapping $F$ is defined everywhere because $p$ has rational components and, hence, there is always $x \in X$ such that $p = x/s(x)$. The mapping $F$ is well defined. Indeed, suppose now there are $x,y$ such that $p = x/s(x)$ and $p = y/s(y)$. By \Hom\ of degree 0, $f(x, \pi) = f(y, \pi)$. Therefore, $F(p)=g(x)=g(y)$.
\end{proof}

We say that a set $S$ in $\QQ^{k}$ is rational convex if whenever $u,v \in S$, then $\alpha u + (1-\alpha) v \in S$ for all rational $\alpha \in [0,1]$.

\begin{lemma}
\label{lemma:paraboloid} 
Let $S$ be a rational convex subset of $\QQ^{2}$ such that  $S$ is full-dimensional. 
Let $g: S \rightarrow \RR_{+}$ be a mapping such that the graph of $g$ is a parabola on any line segment $r \subset S$. Then  $g(u,v) = \rho u^{2} + \sigma v^{2} + \tau uv + \mu u + \nu v+ \xi$ for some real $\rho, \sigma, \tau, \mu, \nu,  \xi$.
\end{lemma}

\begin{proof}
Since $S$ is full-dimensional, the interior of $S$ is not empty and we can suppose without loss of generality that $(0,0) \in \Int S$.
Let us consider the line  defined by $(\alpha t, (1-\alpha) t)$ for some $\alpha \in \QQ$ and all $t \in \QQ$. The intersection of this line with $S$ defines a line segment $r_{\alpha}$ passing by the origin.
The graph of $g$ on $r_{\alpha}$ is a parabola. We can express this by means of the following polynomial of degree 2 in $t$\,:
\begin{equation}
\label{eq:paraboleRayonNew}
g(\alpha t, (1-\alpha) t)  =  k_{\alpha} t^{2}  + l_{\alpha} t + m_{\alpha}, 
\end{equation}
where $k_{\alpha}, l_{\alpha}$ and $m_{\alpha}$ are real numbers.

Let us now consider the line defined by $(\alpha t, (1-\alpha) t)$ for some $t \in \QQ$ and all $\alpha \in \QQ$. The intersection of this line with $S$ defines a line segment $s_{t}$. We can express that the graph of $g$ on $s_{t}$ is a parabola by means of a polynomial of degree 2 in $\alpha$\,:
\begin{equation}
\label{eq:paraboleSecanteNew}
g(\alpha t, (1-\alpha) t) = \alpha^{2} \beta_{t} + \alpha \gamma_{t} + \delta_{t},
\end{equation}
where $\beta_{t}, \gamma_{t}$ and $\delta_{t}$ are real numbers. 
Setting $t=0$ in \eqref{eq:paraboleSecanteNew} yields, $g(0,0)= \alpha^{2} \beta_{0} + \alpha \gamma_{0}+ \delta_{0}$. Since this must be true for all $\alpha$, we must have $\beta_{0} = \gamma_{0} = 0$.

Equating  \eqref{eq:paraboleRayonNew} and \eqref{eq:paraboleSecanteNew} yields 
\begin{equation}
\label{eq:paraboleRayonSecanteNew}
k_{\alpha} t^{2}  + l_{\alpha} t + m_{\alpha} = \alpha^{2} \beta_{t} + \alpha \gamma_{t} + \delta_{t}.
\end{equation}
Setting $t=0$, $t=1$ and $t=2$ in \eqref{eq:paraboleRayonSecanteNew} yields
\begin{eqnarray*}
m_{\alpha} & = &    \delta_{0}   \\
k_{\alpha}  + l_{\alpha}  + m_{\alpha} & = &   \alpha^{2} \beta_{1} + \alpha \gamma_{1} + \delta_{1}   \\
4 k_{\alpha}  + 2 l_{\alpha}  + m_{\alpha} & = &   \alpha^{2} \beta_{2} + \alpha \gamma_{2} + \delta_{2}.   
\end{eqnarray*}
The solution of this system is
\begin{eqnarray*}
m_{\alpha} & = &  \delta_{0}  \\
k_{\alpha}   & = &   \alpha^{2} \ \frac{ \beta_{2} - 2 \beta_{1} }{2} + \alpha \frac{ \gamma_{2} - 2 \gamma_{1} }{2} + \frac{ \delta_{2} - 2 \delta_{1} +  \delta_{0} }{2}  \\
 l_{\alpha}   & = &  \alpha^{2} \ \frac{4 \beta_{1} - \beta_{2} }{2} + \alpha \ \frac{4 \gamma_{1} - \gamma_{2} }{2} + \frac{4 \delta_{1}- \delta_{2} -3  \delta_{0} }{2}.   
\end{eqnarray*}
Let us rewrite \eqref{eq:paraboleRayonNew}\,: 
\begin{eqnarray*}
g(\alpha t, (1-\alpha) t) 
& = &  \left(\alpha^{2} \ \frac{ \beta_{2} - 2 \beta_{1} }{2} + \alpha \ \frac{ \gamma_{2} - 2 \gamma_{1} }{2} + \frac{ \delta_{2} - 2 \delta_{1} +  \delta_{0} }{2} \right) t^{2}  \\
 &  & + \left( \alpha^{2} \ \frac{4 \beta_{1} - \beta_{2} }{2} + \alpha \ \frac{4 \gamma_{1} - \gamma_{2} }{2} + \frac{4 \delta_{1}- \delta_{2} -3  \delta_{0} }{2} \right) t + \delta_{0}.  
\end{eqnarray*}
Letting $\alpha t = u, (1-\alpha) t = v, \alpha = u/(u+v)$ and $t = u+v$, we find that $g(u,v)$ is equal to 
\begin{eqnarray}
\label{eq:g(u,v)}
&  &  u^{2} \ \frac{ \beta_{2} - 2 \beta_{1} }{2} + u(u+v) \ \frac{ \gamma_{2} - 2 \gamma_{1} }{2} + (u+v)^{2} \ \frac{ \delta_{2} - 2 \delta_{1} +  \delta_{0} }{2}   \nonumber \\
 &  & +  \frac{u^{2}}{u+v} \ \frac{4 \beta_{1} - \beta_{2} }{2} + u \ \frac{4 \gamma_{1} - \gamma_{2} }{2} + (u+v) \ \frac{4 \delta_{1}- \delta_{2} -3  \delta_{0} }{2}  + \delta_{0}. \ \ \ \ \ \
 \end{eqnarray}

The graph of $g(u,v)$ must be a parabola on the line segment corresponding to  $v=u+1$. That is,
\begin{multline*}
 u^{2} \ \frac{ \beta_{2} - 2 \beta_{1} }{2} + u(2u+1) \ \frac{ \gamma_{2} - 2 \gamma_{1} }{2} + (2u+1)^{2} \ \frac{ \delta_{2} - 2 \delta_{1} +  \delta_{0} }{2}   \\
  +  \frac{u^{2}}{2u+1} \ \frac{4 \beta_{1} - \beta_{2} }{2} + u \ \frac{4 \gamma_{1} - \gamma_{2} }{2} + (2u+1) \ \frac{4 \delta_{1}- \delta_{2} -3  \delta_{0} }{2}  + \delta_{0}
\end{multline*}
must be a parabola in $u$. This is possible only if $4 \beta_{1} - \beta_{2} = 0$. We have therefore reached the conclusion that \eqref{eq:g(u,v)} can be written as in the statement of the lemma.
\end{proof}

Let $K = \{ 1, 2, \ldots, k \}$ and $K^{*} = \{ 1, 2, \ldots, k-1 \}$.

\begin{lemma}
\label{lemma:paraboloid5}
Let $S$ be a rational convex subset of $\QQ^{k}$ such that  $S$ is full-dimensional.
Let $g: S \rightarrow \RR_{+}$ be a mapping such that  the graph of $g$ is a parabola on any line segment $r \subset S$. 
Suppose the restriction of $g$ to the hyperplane defined by $\sum_{i \in K} u_{i} = t$ ( for all $t \in \RR$ such that the hyperplane intersects $S$) has the form $g(u_{1}, \ldots, u_{k-1}, t - \sum_{i \in K^{*}} u_{i}) = \sum_{i \in K^{*}} \sigma_{ii}u_{i}^{2} + \sum_{i,j \in K^{*}, i < j} \sigma_{ij} u_{i} u_{j} + \sum_{i \in K^{*}} \sigma_{i} u_{i} + \sigma_{0}$ for some real $\sigma_{ii}, \sigma_{ij}, \sigma_{i}, \sigma_{0}$.

Then $g(u_{1}, \ldots, u_{k}) = \sum_{i \in K} \rho_{ii}u_{i}^{2} + \sum_{i,j \in K, i < j} \rho_{ij} u_{i} u_{j} + \sum_{i \in K} \rho_{i} u_{i} + \rho_{0}$ for some real $\rho_{ii}, \rho_{ij}, \rho_{i}, \rho_{0}$.
\end{lemma}

\begin{proof} 
Since $S$ is full-dimensional, there is $u \in \Int S$ and we can suppose without loss of generality that $u=(0, \ldots, 0)$.
Let us consider the line  defined by $(\alpha_{1} t, \alpha_{2} t, \ldots, \alpha_{k-1} t,  (1- \sum_{i \in K^{*}}\alpha_{i}) t)$ for some $\alpha_{1}, \ldots, \alpha_{k-1} \in \RR$ and all $t \in \RR$. The intersection of this line with $S$ defines a line segment $r_{\alpha}$ passing by the origin.
The graph of $g$ on $r_{\alpha}$ is a parabola. We can express this by means of the following polynomial of degree 2 in $t$\,:
\begin{equation}
\label{eq:paraboleRayon5}
g(\alpha_{1} t, \alpha_{2} t, \ldots, \alpha_{k-1} t,  (1- \sum_{i \in K^{*}}\alpha_{i}) t)  =  k_{\alpha } t^{2}  + l_{\alpha } t + m_{\alpha }, 
\end{equation}
where $k_{\alpha}, l_{\alpha}$ and $m_{\alpha}$ are real numbers.

Let us now consider the hyperplane defined by $(\alpha_{1} t, \alpha_{2} t, \ldots, \alpha_{k-1} t,  (1- \sum_{i \in K^{*}}\alpha_{i}) t)$ for some $t \in \RR$ and  $\alpha_{i}  \in \RR, \forall i \in K^{*}$.  We assumed in the statement of the lemma,
\begin{multline}
\label{eq:paraboleSecante5}
g(\alpha_{1} t, \alpha_{2} t, \ldots, \alpha_{k-1} t,  (1- \sum_{i \in K^{*}}\alpha_{i}) t) =   \\
 \sum_{i \in K^{*}} \sigma^{t}_{ii} \alpha_{i}^{2} + \sum_{i,j \in K^{*}, i < j} \sigma^{t}_{ij} \alpha_{i} \alpha_{j} + \sum_{i \in K^{*}} \sigma^{t}_{i} \alpha_{i} + \sigma^{t}_{0}.
\end{multline}
Setting $t=0$ in \eqref{eq:paraboleSecante5} yields, $g(0,\ldots,0)= \sum_{i \in K^{*}} \sigma^{0}_{ii} \alpha_{i}^{2} + \sum_{i,j \in K^{*}, i < j} \sigma^{0}_{ij} \alpha_{i} \alpha_{j} + \sum_{i \in K^{*}} \sigma^{0}_{i} \alpha_{i} + \sigma^{0}_{0}$. Since this must be true for all $\alpha_{i} \in \RR, i \in K^{*}$, we must have $\sigma^{0}_{ii} = \sigma^{0}_{ij} =  \sigma^{0}_{i} = 0$, for all $i,j \in K^{*}$.

Equating  \eqref{eq:paraboleRayon5} and \eqref{eq:paraboleSecante5} yields 
\begin{equation}
\label{eq:paraboleRayonSecante5}
k_{\alpha } t^{2}  + l_{\alpha } t + m_{\alpha } = \sum_{i \in K^{*}} \sigma^{t}_{ii} \alpha_{i}^{2} + \sum_{i,j \in K^{*}, i < j} \sigma^{t}_{ij} \alpha_{i} \alpha_{j} + \sum_{i \in K^{*}} \sigma^{t}_{i} \alpha_{i} + \sigma^{t}_{0}.
\end{equation}
Setting $t=0$, $t=1$ and $t=2$ in \eqref{eq:paraboleRayonSecante5} yields
\begin{eqnarray*}
m_{\alpha } & = &    \sigma_{0}   \\
k_{\alpha }  + l_{\alpha }  + m_{\alpha } & = &   \sum_{i \in K^{*}} \sigma^{1}_{ii} \alpha_{i}^{2} + \sum_{i,j \in K^{*}, i < j} \sigma^{1}_{ij} \alpha_{i} \alpha_{j} + \sum_{i \in K^{*}} \sigma^{1}_{i} \alpha_{i} + \sigma^{1}_{0}   \\
4 k_{\alpha }  + 2 l_{\alpha }  + m_{\alpha } & = &  \sum_{i \in K^{*}} \sigma^{2}_{ii} \alpha_{i}^{2} + \sum_{i,j \in K^{*}, i < j} \sigma^{2}_{ij} \alpha_{i} \alpha_{j} + \sum_{i \in K^{*}} \sigma^{2}_{i} \alpha_{i} + \sigma^{2}_{0}.   
\end{eqnarray*}  
The solution of this system is
\begin{eqnarray*}
m_{\alpha } & = &  \sigma_{0}  \\
k_{\alpha }   & = &  \sum_{i \in K^{*}} \alpha_{i}^{2} \ \frac{ \sigma^{2}_{ii} - 2 \sigma^{1}_{ii}}{2} 
	+ \sum_{i,j \in K^{*}, i < j} \alpha_{i} \alpha_{j} \ \frac{ \sigma^{2}_{ij} - 2 \sigma^{1}_{ij} }{2} \\
 & &	+ \sum_{i \in K^{*}} \alpha_{i} \  \frac{ \sigma^{2}_{i} - 2 \sigma^{1}_{i} }{2} +  \frac{ \sigma^{2}_{0} - 2 \sigma^{1}_{0} + \sigma^{0}_{0} }{2}  \\
 l_{\alpha }   & = &   \sum_{i \in K^{*}} \alpha_{i}^{2} \ \frac{ 4 \sigma^{1}_{ii} -  \sigma^{2}_{ii}}{2} 
	+ \sum_{i,j \in K^{*}, i < j} \alpha_{i} \alpha_{j} \ \frac{ 4 \sigma^{1}_{ij} -  \sigma^{2}_{ij} }{2} \\
 & &	+ \sum_{i \in K^{*}} \alpha_{i} \  \frac{ 4 \sigma^{1}_{i} -  \sigma^{2}_{i}  }{2} +  \frac{ 4 \sigma^{1}_{0} -  \sigma^{2}_{0} - 3 \sigma^{0}_{0} }{2} .   
\end{eqnarray*}
Let us rewrite \eqref{eq:paraboleRayon5}\,: $g(\alpha_{1} t, \alpha_{2} t, \ldots, \alpha_{k-1} t,  (1- \sum_{i \in K^{*}}\alpha_{i}) t) = $
\begin{align*}
  \left(  \sum_{i \in K^{*}} \alpha_{i}^{2} \ \frac{ \sigma^{2}_{ii} - 2 \sigma^{1}_{ii}}{2} \right.
   &+ \sum_{i,j \in K^{*}, i < j} \alpha_{i} \alpha_{j} \ \frac{ \sigma^{2}_{ij} - 2 \sigma^{1}_{ij} }{2}  
   \\
 	& \left. + \sum_{i \in K^{*}} \alpha_{i} \  \frac{ \sigma^{2}_{i} - 2 \sigma^{1}_{i} }{2} +  \frac{ \sigma^{2}_{0} - 2 \sigma^{1}_{0} + \sigma^{0}_{0} }{2}  \right) t^{2}  \\
  + \left( \sum_{i \in K^{*}} \alpha_{i}^{2} \ \frac{ 4 \sigma^{1}_{ii} -  \sigma^{2}_{ii}}{2} \right.
	&+ \sum_{i,j \in K^{*}, i < j} \alpha_{i} \alpha_{j} \ \frac{ 4 \sigma^{1}_{ij} -  \sigma^{2}_{ij} }{2} 
 	  \\
	& \left. + \sum_{i \in K^{*}} \alpha_{i} \  \frac{ 4 \sigma^{1}_{i} -  \sigma^{2}_{i}  }{2} 
	+  \frac{ 4 \sigma^{1}_{0} -  \sigma^{2}_{0} - 3 \sigma^{0}_{0} }{2}   \right) t + \sigma_{0}.  
\end{align*}
Letting $\alpha_{i} t = u_{i}, \forall i \in K^{*}, (1-\sum_{i \in K^{*}} \alpha_{i}) t = u_{k}$, we have $\alpha_{i} = u_{i}/\sum_{i \in K} u_{i}$ and $t = \sum_{i \in K} u_{i}$, and the previous equation becomes, $g(u_{1}, \ldots, u_{k}) =$
\begin{multline*}
 \sum_{i \in K^{*}} u_{i}^{2} \ \frac{ \sigma_{ii}^{2} - 2 \sigma_{ii}^{1} }{2} 
+ \sum_{i,j \in K^{*}, i < j} u_{i} u_{j} \ \frac{ \sigma^{2}_{ij} - 2 \sigma^{1}_{ij} }{2} 
+ \sum_{i \in K^{*}} u_{i} \sum_{j \in K} u_{j} \ \frac{  \sigma^{2}_{i} - 2 \sigma^{1}_{i} }{2} \\
   +  \left(\sum_{i \in K} u_{i} \right)^{2}   \ \frac{\sigma^{2}_{0} - 2 \sigma^{1}_{0} + \sigma^{0}_{0} }{2}   
  + \sum_{i \in K^{*}}  \frac{u_{i}^{2}}{\sum_{j \in K} u_{j}} \ \frac{4 \sigma^{1}_{ii} -  \sigma^{2}_{ii} }{2}  
   \\
 +   \sum_{i,j \in K^{*}, i < j} \frac{u_{i} u_{j}}{\sum_{j \in K} u_{j}} \ \frac{4 \sigma^{1}_{ij} -  \sigma^{2}_{ij}  }{2}  
 + \sum_{i \in K^{*}}  u_{i} \ \frac{4 \sigma^{1}_{i} -  \sigma^{2}_{i} }{2}  \\
  + \sum_{i \in K} u_{i} \ \frac{4 \sigma^{1}_{0} -  \sigma^{2}_{0} - 3 \sigma^{0}_{0}}{2}  + \sigma_{0}.
 \end{multline*}

For any $j \in K^{*}$, the graph of $g(u_{1}, \ldots, u_{k})$ must be a parabola on the line segment corresponding to  $u_{i}=0, \forall i \in K^{*}\setminus \{j\} , u_{k}=u_{j}+1$. That is, $g(0, \ldots, 0, u_{j},0, \ldots, 0,u_{j}+1) =$
\begin{multline*}
  u_{j}^{2} \ \frac{ \sigma_{jj}^{2} - 2 \sigma_{jj}^{1} }{2} 
+  u_{j} (2 u_{j} + 1) \ \frac{  \sigma^{2}_{j} - 2 \sigma^{1}_{j} }{2} 
+ \left( 2 u_{j} + 1 \right)^{2} \ \frac{\sigma^{2}_{0} - 2 \sigma^{1}_{0} + \sigma^{0}_{0} }{2}  \\    
  +   \frac{u_{j}^{2}}{2 u_{j} + 1} \ \frac{4 \sigma^{1}_{jj} -  \sigma^{2}_{jj} }{2}  
 +  u_{j} \ \frac{4 \sigma^{1}_{j} -  \sigma^{2}_{j} }{2} 
  + (2 u_{j} + 1) \ \frac{4 \sigma^{1}_{0} -  \sigma^{2}_{0} - 3 \sigma^{0}_{0}}{2}  + \sigma_{0}
 \end{multline*}
must be a parabola in $u_{j}$. This is possible only if $4 \sigma^{1}_{jj} -  \sigma^{2}_{jj} = 0$ for all $j \in K^{*}$. 

Similarly, for any $i,j \in K^{*}$ with $i < j$, the graph of $g(u_{1}, \ldots, u_{k})$ must be a parabola on the line segment corresponding to  $u_{i}=u_{j}$, $u_{l}=0, \forall l \in K^{*} \setminus\{i,j\} , u_{k}=u_{i}+1$. That is, $g(0, \ldots, 0,  u_{i}, 0, \ldots, 0, u_{i}, 0, \ldots, 0, u_{i}+1) =$
\begin{multline*}
  u_{i}^{2} \ \frac{ \sigma_{ii}^{2} - 2 \sigma_{ii}^{1} }{2}  + u_{i}^{2} \ \frac{ \sigma_{jj}^{2} - 2 \sigma_{jj}^{1} }{2} 
+ u_{i}^{2} \ \frac{ \sigma^{2}_{ij} - 2 \sigma^{1}_{ij} }{2} 
+ u_{i} (3 u_{i} + 1) \ \frac{  \sigma^{2}_{i} - 2 \sigma^{1}_{i} }{2}  \\
+ u_{i} (3 u_{i} + 1) \ \frac{  \sigma^{2}_{j} - 2 \sigma^{1}_{j} }{2}
   + \left(  3 u_{i} + 1 \right)^{2} \ \frac{\sigma^{2}_{0} - 2 \sigma^{1}_{0} + \sigma^{0}_{0} }{2}   
   +   \frac{u_{i}^{2} }{3 u_{i} + 1} \ \frac{4 \sigma^{1}_{ij} -  \sigma^{2}_{ij}  }{2}   \\
 +   u_{i} \ \frac{4 \sigma^{1}_{i} -  \sigma^{2}_{i} }{2} 
 +u_{i} \ \frac{4 \sigma^{1}_{j} -  \sigma^{2}_{j} }{2} 
  + (3 u_{i} + 1) \ \frac{4 \sigma^{1}_{0} -  \sigma^{2}_{0} - 3 \sigma^{0}_{0}}{2}  + \sigma_{0}
 \end{multline*}
 must be a parabola in $u_{i}$. This is possible only if $4 \sigma^{1}_{ij} -  \sigma^{2}_{ij} = 0$ for all $i,j \in K^{*}$ with $i \neq j$.
We have therefore reached the conclusion that $g$ has the desired form.
\end{proof}

\begin{lemma}
\label{lemma:paraboloid6}
Let $S$ be a rational convex subset of $\QQ^{k}$ such that  $(0,0, \ldots, 0) \in \Int S$.
Let $g: S \rightarrow \RR_{+}$ be a mapping such that $g(u_{1}, \ldots, u_{k}) = 0$ if and only  $(u_{1}, \ldots, u_{k})=(0, \ldots,0)$ and the graph of $g$ is a parabola on any line segment $r \subset S$. Then $g(u_{1}, \ldots, u_{k}) = \sum_{i \in K} \rho_{ii}u_{i}^{2} + \sum_{i,j \in K, i < j} \rho_{ij} u_{i} u_{j}$ for some real $\rho_{ii}, \rho_{ij}$.
\end{lemma}

\begin{proof}
By induction and  Lemmas~\ref{lemma:paraboloid5} and ~\ref{lemma:paraboloid} , $g(u_{1}, \dots, u_{k}) = \sum_{i \in K} \rho_{ii}u_{i}^{2} + \sum_{i,j \in K, i < j} \rho_{ij} u_{i} u_{j} + \sum_{i \in K} \rho_{i} u_{i} + \rho_{0}$ for some real $\rho_{ii}, \rho_{ij}, \rho_{i}, \rho_{0}$.
On the line $u_{2}=u_{3} = \ldots = u_{k}=0$, $g(u_{1},0, \ldots, 0) = \rho_{11}u^{2} + \rho_{1} u + \rho_{0}$. The graph of this function of $u$ must be a parabola with vertex in 0. Hence $\rho_{1} = \rho_{0}=0$. Similarly, for any $i \in K$, considering the line defined by $u_{j}=0, \forall j \neq i$  entails $\rho_{i}=0$. In conclusion, $g(u_{1}, \ldots, u_{k}) = \sum_{i \in K} \rho_{ii}u_{i}^{2} + \sum_{i,j \in K, i < j} \rho_{ij} u_{i} u_{j}$. 
\end{proof}

\begin{lemma}
Suppose $\pi$ is fixed and the dissimilarity measure $f$ satisfies  \Hom\ of degree 0 and \AQ.   Then, for all $p, q \in \Pi$ such that  $p + q - \pi \in \Pi$ and $p - q + \pi \in \Pi$, we have
\begin{equation}
\label{Q'Add1}
F(p + q - \pi) + F(p - q + \pi) = 2 F(p) + 2 F(q).
\end{equation} 
\end{lemma}

\begin{proof}
Let $p,q$ be as in the statement of the lemma.
There are two distributions $x,y \in X$ such that $x/s(x) = p$ and $y/s(y) = q$. Hence $F(p)=f(x,\pi)$ and $F(q)=f(y,\pi)$. By  \Hom\ of degree 0, $F(p)=f(x,\pi) = f(s(y)x,\pi)$ and $F(q)=f(y,\pi) = f(s(x)y,\pi)$. The two distributions $s(y)x$ and $s(x)y$ have the same size, i.e., $s(x)s(y)$. 
Hence we can apply \AQ\ and we find 
\begin{eqnarray*}
2 F(p) + 2 F(q) & = & 2 f(s(y)x, \pi) + 2 f(s(x)y,\pi) \\
                        & = & f(s(y)x + s(x)y - s(x)s(y) \pi, \pi)  \\
                        &    & + f(s(y)x - s(x)y + s(x)s(y) \pi, \pi).
\end{eqnarray*}
By the definition of $F$, 
$$f(s(y)x + s(x)y - s(x)s(y) \pi, \pi) = F(p+q-\pi)$$
and
$$f(s(y)x - s(x)y + s(x)s(y) \pi, \pi) =F(p-q+\pi).$$
In conclusion,
$$F(p + q - \pi) + F(p - q + \pi) = 2 F(p) + 2 F(q). \eqno \qedhere$$
\end{proof}

For every $l \in N$, define $N_{l} = N \setminus \{l\}$ and $N_{lm} = N \setminus \{l, m\}$.

\begin{lemma}
\label{lemme.chi2. paraboloid}
Suppose $\pi$ is fixed and the dissimilarity measure $f$ satisfies  \Hom\ of degree 0 and \AQ.    Then, for every $l \in N$ and $p \in \Pi$,  
$$F(p) =   \sum_{i \in N_{l}} \rho^{l}_{ii}(p_{i}-\pi_{i})^{2} + \sum_{i,j \in N_{l} : i < j} \rho^{l}_{ij} (p_{i} - \pi_{i}) (p_{j} - \pi_{j})$$ for some real $\rho^{l}_{ii}, \rho^{l}_{ij}.$
\end{lemma}

\begin{proof}
Let $r$ be a line segment with extremities $s,t \in \Pi$, with $s \neq t$.
Every point of $r \cap \Pi$ can be written as $\alpha t+ (1-\alpha) s$ with $\alpha \in [0,1]$ and $\alpha$ rational. 
Consider any two points $p,q \in r \cap \Pi$, the position of which on $r$ is characterized by $\alpha$ and $\beta$ respectively. Then $p + q - \pi$ lies on the line segment between $2t - \pi$ and $2s - \pi$ and it can be written as $$\left( \frac{\alpha + \beta}{2} \right) (2t - \pi) + \left( 1-\frac{\alpha + \beta}{2} \right) (2s - \pi).$$
Similarly, $p - q + \pi$ lies on the line segment between $t-s + \pi$ and $s - t + \pi$ and it can be written as $$\left( \frac{\alpha - \beta + 1}{2} \right) (t -s + \pi) + \left( 1-\frac{\alpha - \beta + 1}{2} \right) (s - t + \pi).$$
Notice that $(\alpha+\beta)/2 \in [0,1]$ and $(\alpha-\beta+1)/2 \in [0,1]$ for any $\alpha, \beta \in [0,1]$.
Define three mappings as follows\,:
\begin{itemize}

  \item     $L : [0,1] \cap \mathbb{Q} \rightarrow \mathbb{R}_{+}$ by $L(\alpha) = F(p)$ if $p = \alpha t + (1-\alpha) s$;
  
  \item    $G : [0,1] \cap \mathbb{Q} \rightarrow \mathbb{R}_{+}$ by $G(\alpha) = F(p)$ if $p = \alpha (2t - \pi) + (1-\alpha) (2s - \pi)$;
  
  \item   $H : [0,1] \cap \mathbb{Q} \rightarrow \mathbb{R}_{+}$ by $H(\alpha) = F(p)$ if $p = \alpha (t-s + \pi) + (1-\alpha) (s-t + \pi)$.
  
\end{itemize}
Then \eqref{Q'Add1} can be rewritten as 
\begin{equation}
\label{eq.fonctionelle.paraboloide}
G\left( \frac{\alpha + \beta}{2}  \right) + H\left( \frac{\alpha - \beta + 1}{2}  \right) = 2 L(\alpha) + 2 L(\beta)
\end{equation}
and it holds for all rational $\alpha, \beta \in [0,1]$. This functional equation is a generalization of Equation (18) discussed in \cite[p.82]{Aczel66}.

If $\alpha = \beta$, then $G(\alpha) + H(1/2) = 4 L(\alpha)$. In other words, $G(\alpha) = 4 L(\alpha) + \delta'$ for some real number $\delta'$. If $\alpha = 1 - \beta$, then $G(1/2)+H(1-\beta) = 2 L(1-\beta) + 2 L(\beta)$. So, $H(1-\beta) = 2 L (\beta) + 2 L (1-\beta) + \delta''$ for some real number $\delta''$. We can now rewrite \eqref{eq.fonctionelle.paraboloide} as
\begin{equation*}
\label{eq.fonctionelle.paraboloide2}
4 L\left( \frac{\alpha + \beta}{2}  \right) + 2 L \left( \frac{\alpha - \beta + 1}{2}  \right) + 2 L\left( \frac{\beta - \alpha + 1}{2}  \right)  = 2 L(\alpha) + 2 L(\beta) + \delta'''
\end{equation*}
with $\delta''' = -\delta' - \delta''$. If we now let $\alpha = (m-2)c$ and $\beta = mc$ with $m$ a positive integer ($m \geq 2$) and $c$ a positive rational number such that $mc \in [0,1]$, then 
$$4 L\left( (m-1)c  \right) + 2 L \left( (1/2)-c  \right) + 2 L \left( (1/2)+c \right)  = 2 L ((m-2)c) + 2 L(mc) + \delta'''.$$
If we divide this equation by 2 and reorder the terms (with $\delta = \delta'''/2$), we obtain that 
$L(mc) - L((m-1)c)$ is equal to 
\begin{eqnarray*}
 &   & L((m-1)c) - L((m-2)c) + L((1/2)-c) + L((1/2)+c) + \delta \\
 & = & L((m-2)c) - L((m-3)c) + 2 \big(L((1/2)-c) + L((1/2)+c) + \delta \big) \\
 & = & L((m-3)c) - L((m-4)c) + 3 \big(L((1/2)-c) + L((1/2)+c) + \delta \big) \\
 & = & \ldots \\
 & = & L(c) - L(0) + (m-1) \big(L((1/2)-c) + L((1/2)+c) + \delta \big).
\end{eqnarray*}
Notice that, for all $m \in \NN (m \geq 2)$ and $c \in \QQ_{++}$ such that $mc \in [0,1]$,
\begin{eqnarray}
L(mc) & = & \sum_{i=1}^{m} \big( L(ic) - L((i-1)c) \big) + L(0)   \nonumber \\
          & = & \sum_{i=1}^{m} \Big(L(c) - L(0) + (i-1) \big(L((1/2)-c) + L((1/2)+c) + \delta \big) \Big) \nonumber   \\
          &    & \ \ \ \ \ + L(0)  \nonumber   \\
          & = & m L(c) + (1-m) L(0) \nonumber \\
          &    & \ \ \ \ \ + \frac{m(m-1)}{2} \ \big(L((1/2)-c) + L\big((1/2)+c ) + \delta \big).\label{eq:parabola}
\end{eqnarray}
If $m=0$, it is easy to verify that \eqref{eq:parabola} still holds. Indeed, 
$$L(0c)= 0L(c)+(1-0)L(0)+\frac{0(0-1)}{2} \ \big(L((1/2)-c) + L\big((1/2)+c ) + \delta \big).$$
A similar reasoning holds when $m=1$.
Hence, for any fixed value of $c$ (a positive rational number), $L(\beta)$ is a polynomial  of degree 2 in $\beta$ for every $\beta \in [0,1]$ that can be written as $mc$ for some integer $m$, hence we can also write \eqref{eq:parabola} as 
$$
L(\beta)=  a \beta^{2} + b \beta + d,
$$
for some real numbers $a, b, d$.
Let $c=\sfrac{1}{2}$; $\beta$ can take the values $0, \sfrac{1}{2}, 1$. We have 
$$
\left\{
\begin{array}{lcl}
  L(0) & = & d \\
  L(\sfrac{1}{2}) & = & \sfrac{1}{4}\, a + \sfrac{1}{2}\, b + d \\
  L(1) & = & a+ b+ d 
\end{array}
\right.
$$
This is a non-singular (determinant = $-1/4$) system of linear equations which determines unique values for $a, b, d$. These are a linear combination of the values of $L(\beta)$ for $\beta = 0, \sfrac{1}{2}, 1$. There is a single parabola that passes through the three points $(\beta, L(\beta))$ for $\beta = 0, \sfrac{1}{2}, 1$.

Consider now any rational number $\beta=\frac{w}{w'}$. If $w'$ is odd, we also have that $\beta = \frac{2w}{2w'}$ so that we can assume that $w'$ is even. Using \eqref{eq:parabola}, we have that $L(\beta)=  a' \beta^{2} + b' \beta + d'$ for some $a', b', d'$ and for any integer $w$ in the interval $[0,w']$. We have to show that these constants are $a, b, d$. Indeed, for $w=0, w'/2, w'$, we have that $\beta = 0, \sfrac{1}{2}, 1$ respectively. Hence, the points $(\beta, L(\beta))$ for $\beta = 0, \sfrac{1}{2}, 1$ are also on the parabola with coefficients $a', b', c'$ and these points are the same as in the case $w'=2$ since $L$ is a single function. Since there is only one parabola through these points and the parabola with coefficients $a, b , c$ passes through these points, we conclude that $a=a'$, $b=b'$ and $d=d'$.

This being true for any (even) value of $w'$, we have that $L(\beta)$ is a quadratic function of $\beta$ independently of the denominator of the rational number $\beta$. 

Of course, the coefficients of the polynomial $a,b$ and $d$ depend on the line segment $r$ joining $t$ and $s$. So, we had better write $L_{r}(\beta) = a_{r} \beta^{2} + b_{r} \beta + d_{r}$. We now go back to $F$. Since $L_{r}$ is a polynomial of degree 2 for any $r$, we find that the graph of $F$ is a parabola on any line segment $r$. Define $G: \Pi \rightarrow \RR_{+}$ by $G(p-\pi) = F(p)$. Then $G$ is like $g$ in the statement of Lemma~\ref{lemma:paraboloid6}, with $k=n-1$ (because $\Pi$ has only $n-1$ dimensions).
  Then, $G(p-\pi) = F(p) = \sum_{i \in N_{l}} \rho^{l}_{ii}(p_{i}-\pi_{i})^{2} + \sum_{i,j \in N_{l}, i < j} \rho^{l}_{ij} (p_{i} - \pi_{i}) (p_{j} - \pi_{j})$ for some real $\rho^{l}_{ii}, \rho^{l}_{ij}$. Notice that, for each $l \in N$, we have such an expression.   
\end{proof}

\bigskip
\noindent {\bf Proof of Theorem~\ref{Chi2}.} \hspace{0.2 em} 
\paragraph{Case $n=2$.}
\subparagraph{$f$ is homogeneous of degree 0.}
Since $n=2$, we have $(p_{1}-\pi_{1})^{2} = (p_{2}-\pi_{2})^{2}$. Then Lemma~\ref{lemme.chi2. paraboloid} yields $F(p)=\rho_{11}^{2} (p_{1}-\pi_{1})^{2} = \frac{\rho_{11}^{2}}{2} \big( (p_{1}-\pi_{1})^{2} + (p_{2}-\pi_{2})^{2} \big)$ and
\begin{equation}
\label{chi.n=2}
f(x,\pi) =  \frac{\rho_{11}^{2}}{2} \left( \left( \frac{x_{1}}{s(x)}-\pi_{1} \right)^{2} + \left( \frac{x_{2}}{s(x)}-\pi_{2} \right)^{2} \right).
\end{equation}

\subparagraph{$f$ is homogeneous of degree 1.}
By Lemma~\ref{lemma:0-1-Homogeneity}, $f/s(x)$ satisfies \Hom\ of degree 0 and \AQ. By \eqref{chi.n=2}, $f/s(x)$ is proportional to $\chi^{2}_{0}$ and $f$  to $\chi^{2}_{1}$.

\paragraph{Case $n \geq 3$.}
\subparagraph{$f$ is homogeneous of degree 0.}
 We know from the previous lemma that $F$ can be expressed as 
\begin{equation}
\label{eq:F-avec-rho}
F(p) =   \sum_{i \in N_{n}} \rho^{n}_{ii}(p_{i}-\pi_{i})^{2} + \sum_{i,j \in N_{n} : i < j} \rho^{n}_{ij} (p_{i} - \pi_{i}) (p_{j} - \pi_{j})
\end{equation}
 for some real $\rho^{n}_{ii}, \rho^{n}_{ij}$. Recall that $N_{n}= N \setminus \{n\}$. So, for all $j, l \in N_{n}$, with $j \neq l$,
\begin{equation*}
\label{eq:f(kPi+1-1)}
f(k \pi + \onej - \onel, \pi) =   \frac{\rho_{jj}^{n} +  \rho_{ll}^{n} -  \rho_{jl}^{n}}{k^{2}},
\end{equation*}
and
\begin{equation*}
\label{eq:f(kPi+1n-1)}
f(k \pi + \onen - \onej, \pi) =   \frac{\rho_{jj}^{n} }{k^{2}}.
\end{equation*}
For the sake of simplicity, for all $j, l \in N$, let $A_{jl}$ denote $\frac{\pi_{j} + \pi_{l}}{\pi_{j} \pi_{l}}$.
Then, thanks to  \RIE, we can write, for all  $j,l \in N_{n}$ with $j \neq l$.
$$
\frac{f(k \pi + \onej - \onel, \pi)}{f(k \pi + \oneone - \onetwo, \pi)} = 
\frac{\rho_{jj}^{n} +  \rho_{ll}^{n} -  \rho_{jl}^{n}}{\rho_{11}^{n} +  \rho_{22}^{n} -  \rho_{12}^{n}} = \frac{A_{jl}}{A_{12}},
$$
$$
\frac{f(k \pi + \onej - \onen, \pi)}{f(k \pi + \oneone - \onen, \pi)} = 
\frac{\rho_{jj}^{n} }{\rho_{11}^{n} } = \frac{A_{jn}}{A_{1n}}
$$
and
$$
\frac{f(k \pi + \oneone - \onen, \pi)}{f(k \pi + \oneone - \onetwo, \pi)} = 
\frac{\rho_{11}^{n}}{\rho_{11}^{n} +  \rho_{22}^{n} -  \rho_{12}^{n}} = \frac{A_{1n}}{A_{12}}.
$$

Using each of these three equations separately, we find 
\begin{equation}
\label{eq:jl-12}
\rho_{jl}^{n} = \rho_{jj}^{n} + \rho_{ll}^{n} + \frac{A_{jl}}{A_{12}} (\rho_{12}^{n} - \rho_{11}^{n} - \rho_{22}^{n}), \ \forall j,l \in N_{n} : j \neq l,
\end{equation}
\begin{equation}
\label{eq:jj-11}
\rho_{jj}^{n} = \rho_{11}^{n} \frac{A_{jn}}{A_{1n}}, \ \forall j \in N_{n}
\end{equation}
and
\begin{equation}
\label{eq:11-12}
\frac{\rho_{11}^{n} +  \rho_{22}^{n} -  \rho_{12}^{n}}{A_{12}} = \frac{\rho_{11}^{n}}{A_{1n}}.
\end{equation}
If we substitute \eqref{eq:jj-11} and \eqref{eq:11-12} in \eqref{eq:jl-12}, we obtain
\begin{eqnarray*}
\rho_{jl}^{n} & = & \rho_{11}^{n} \left( \frac{A_{jn}}{A_{1n}} + \frac{A_{ln}}{A_{1n}} - \frac{A_{jl}}{A_{1n}} \right) , \ \forall j,l \in N_{n} : j \neq l, \\
	& = & 2 \rho_{11}^{n} \ \frac{\pi_{1}}{\pi_{1} + \pi_{n}}.
\end{eqnarray*}
From \eqref{eq:jj-11}, we find 
$$
\rho_{jj}^{n} = \rho_{11}^{n} \frac{\pi_{j} + \pi_{n}}{\pi_{j}} \ \frac{\pi_{1}}{\pi_{1} + \pi_{n}} , \ \forall j \in N_{n}.
$$
Let us define 
$$\gamma = \rho_{11}^{n} \ \frac{\pi_{1} \pi_{n}}{\pi_{1} + \pi_{n}}.$$
Then,
$$\rho_{jl}^{n} = 2 \gamma \ \frac{\pi_{1} + \pi_{n}}{\pi_{1} \pi_{n}} \ \frac{\pi_{1}}{\pi_{1} + \pi_{n}}
= 2 \ \frac{\gamma}{\pi_{n}}, \ \forall j,l \in N_{n} : j \neq l,$$
and 
$$\rho_{jj}^{n} = \gamma \ \frac{\pi_{1} + \pi_{n}}{\pi_{1} \pi_{n}} \ \frac{\pi_{j} + \pi_{n}}{\pi_{j}} \ \frac{\pi_{1}}{\pi_{1} + \pi_{n}} =  \gamma \ \frac{\pi_{j} + \pi_{n}}{\pi_{j} \pi_{n}}, \ \forall j \in N_{n}.$$

If we now substitute the expressions of $\rho_{jl}^{n}$ and $\rho_{jj}^{n} $   into \eqref{eq:F-avec-rho}, we obtain
 $f(x,\pi) $
\begin{eqnarray}
&=& \sum_{i \in N_{n}} \left( \frac{\gamma}{\pi_{i}} + \frac{\gamma}{\pi_{n}} \right) \left(\pi_{i} - \frac{x_{i}}{s(x)} \right)^{2} \nonumber \\
 & & + \sum_{i,j \in N_{n} : i < j} \frac{2 \gamma}{\pi_{n}}  \left(\pi_{i} - \frac{x_{i}}{s(x)} \right) \left(\pi_{j} - \frac{x_{j}}{s(x)} \right)  \nonumber \\
&=&  \sum_{i \in N_{n}}  \frac{\gamma}{\pi_{i}}  \left(\pi_{i} - \frac{x_{i}}{s(x)} \right)^{2}
+  \sum_{i \in N_{n}}  \frac{\gamma}{\pi_{n}}  \left(\pi_{i} - \frac{x_{i}}{s(x)} \right)^{2} \nonumber \\
& & +   \sum_{i,j \in N_{n} : i < j} \frac{2 \gamma}{\pi_{n}}  \left(\pi_{i} - \frac{x_{i}}{s(x)} \right) \left(\pi_{j} - \frac{x_{j}}{s(x)} \right) \nonumber  \\
&=&  \sum_{i \in N_{n}}  \frac{\gamma}{\pi_{i}}  \left(\pi_{i} - \frac{x_{i}}{s(x)} \right)^{2}
+ \frac{\gamma}{\pi_{n}}  \left( \sum_{i \in N_{n}}    \left(\pi_{i} - \frac{x_{i}}{s(x)} \right) \right)^{2}  \nonumber \\
&=&  \sum_{i \in N_{n}}  \frac{\gamma}{\pi_{i}}  \left(\pi_{i} - \frac{x_{i}}{s(x)} \right)^{2}
+ \frac{\gamma}{\pi_{n}}  \left( 1 -  \sum_{i \in N_{n}} \pi_{i} -1 + \sum_{i \in N_{n}} \frac{x_{i}}{s(x)} \right)^{2} \nonumber \\
&=&  \sum_{i \in N_{n}}  \frac{\gamma}{\pi_{i}}  \left(\pi_{i} - \frac{x_{i}}{s(x)} \right)^{2}
+ \frac{\gamma}{\pi_{n}}  \left( \pi_{n} - \frac{x_{n}}{s(x)} \right)^{2} \nonumber \\
&=& \sum_{i \in N}  \frac{\gamma}{\pi_{i}}  \left(\pi_{i} - \frac{x_{i}}{s(x)} \right)^{2} = \gamma \chi^{2}_{0}(x, \pi). \label{eq:chi2-H0}
\end{eqnarray}

\subparagraph{$f$ is homogeneous of degree 1.}
By Lemma~\ref{lemma:0-1-Homogeneity}, $f/s(x)$ satisfies \Hom\ of degree 0, \AQ\ and \IE. By \eqref{eq:chi2-H0}, $f/s(x)$ is proportional to $\chi^{2}_{0}$ and $f$ to $\chi^{2}_{1}$.
\epr

\bigskip
\noindent {\bf Proof of Theorem~\ref{Chi2.pi.variable}.} \hspace{0.2 em} 
We prove the result only for $\chi^{2}_{1}$ (the case of $\chi^{2}_{0}$ being similar).
We want to prove that there exists $\gamma >0$ such that $f(x,\pi)= \gamma \chi^{2}_{1}(x,\pi)$ for all $x \in X$ and all $\pi \in \Pi$. From Theorem~\ref{Chi2}, we have $f(x,\pi)= \gamma_{\pi} \chi^{2}_{1}(x,\pi)$, where we now  use the notation $\gamma_{\pi}$ to make clear that $\gamma_{\pi}$ can depend on $\pi$. 
Thanks to \IE, for all $\pi, \pi' \in \Pi$, we have
$$\frac{f(k \pi + \onej -\onel,\pi)}{f(k \pi' +\onej - \onel,\pi')} = \frac{\gamma_{\pi}\left( \frac{1}{\pi_{j}} + \frac{1}{\pi_{l}}\right)}{\gamma_{\pi'}\left(\frac{1}{\pi'_{j}} + \frac{1}{\pi'_{l}} \right)}  = \frac{\frac{1}{\pi_{j}} + \frac{1}{\pi_{l}}}{\frac{1}{\pi'_{j}} + \frac{1}{\pi'_{l}}}.$$
This is possible only if $\gamma_{\pi}=\gamma_{\pi'}$ for all $\pi, \pi' \in \Pi$. 
\epr





\end{document}